\documentclass{article}
\usepackage{amssymb,latexsym,amsmath,amsthm,mathrsfs}
\usepackage{graphicx}
\newcommand{\comment}[1]{}

\DeclareMathOperator{\Tan}{tang.}
\DeclareMathOperator{\Cot}{cot.}
\DeclareMathOperator{\etc}{etc.}
\DeclareMathOperator{\Cos}{cos.}
\DeclareMathOperator{\Sin}{sin.}
\DeclareMathOperator{\Sec}{sec.}
\begin{document}
\title{Annotations to a certain passage of Descartes for finding the quadrature of the circle\footnote{Presented to the
Berlin Academy on July 20, 1758.
Originally published as
{\em Annotationes in locum quendam Cartesii ad circuli quadraturam spectantem}, Novi Commentarii academiae scientiarum
Petropolitanae \textbf{8} (1763), 157--168. Summarium is pp. 24--27.
E275 in the Enestr{\"o}m index.
Translated from the Latin by Jordan Bell, University of Toronto, Toronto,
Ontario, Canada.
Email: jordan.bell@gmail.com}}
\author{Leonhard Euler}
\date{}
\maketitle

\begin{center}
{\Large Summarium}
\end{center}
{\small
That the circumference of a circle is incommensurable with its diameter,
or that no measure can be given that simultaneously measures both the diameter
and the circumference of the circle,\footnote{Translator: cf. book X, def. I.1 of
Euclid.}
has been observed already by the Geometers of antiquity, 
yet still now it cannot be demonstrated any further than that
all attempts at finding a measure of this type have been in vain.
Namely,
no two numbers can be exhibited which hold between themselves
the exact same ratio as occurs between the diameter and the circumference.
Thus in practice one is wont to use numbers that closely approach this ratio,
examples of which are 7 to 22 of Archimedes, and 113 to 355 of
Metius;
indeed this true ratio has been expressed then by others more accurately with
numbers,\footnote{Translator: Here by ``numbers'' Euler means ``integers''; cf. book VII, def. 2 of Euclid.}
so that even in the computation of the largest circles the error is negligible.
Knowing they are incommensurable though does not
in itself prevent the ratio of the diameter to the circumference
from being assigned geometrically, since the diagonal of a square
is also incommensurable to the side, and in general all irrational quantities
which arise from the extraction of roots can be constructed 
geometrically.\footnote{Translator: I do not know if this means
roots or just square roots. Indeed finding square roots can be done
geometrically (with straightedge and compass), but finding the cube root cannot.}
It seems that the true circumference of a circle belongs to
a very high class of irrationalities, which can only be reached by
repeating infinitely the process of extracting a root; here the best that
can be done geometrically is to express more and more closely
the true ratio of the circumference to the
diameter.
And the Cartesian construction that the Celebrated Author deals with here
works in this way, such that by the continual apposition of rectangles
which decrease according to a certain rule, a line is drawn which is finally
equal to the circumference of a
circle.
This construction has been so ingeniously devised that by its simple
use one is led quickly to the truth, and we should admire this extraordinary
monument to the great insight of its discoverer.
Euler, while saving this discovery from oblivion, also
publishes many singular formulas and series pertaining to measuring
the circle, by which geometric approximations of this kind can be applied to a
greater extent and further ones can be found.
For instance,
he has demonstrated that with $q$ denoting the length of a quadrant of a circle
whose radius is $=1$,
\[
q=\Sec\frac{1}{2}q \cdot \Sec\frac{1}{4}q\cdot \Sec\frac{1}{8}q\cdot
\Sec\frac{1}{16}q\cdot \Sec\frac{1}{32}q\cdot \etc,
\]
from which one can conclude
the following rather neat and elegant
construction
With the quadrant $AOB$ established, 
the normal $BC$ to the radius $OB$ intersects the line $OC$ bisecting
the angle $AOB$ at $C$. Then, the normal $CD$ at $C$ to this $OC$ intersects
the line $OD$ bisecting the angle $AOC$ at $D$.
Similarly, $DE$, the line normal to this $OD$, intersects the line
$OE$ bisecting the angle $AOD$ at $E$. Again then,
$EF$, the line normal to $OE$, intersects the line $OF$ bisecting the angle $AOE$
at $F$, and so on. 
One continues in this manner until finally the radius $OA$ is reached;
this construction finally stops at the point $Z$. Then having done this,
the line $OZ$ will be precisely equal to the length of the
quadrant $BcdefgA$.\footnote{Translator: Euler in fact does not explain why $OZ=BcdefgA$
later in the paper. The following explanation is from Ed Sandifer. $OBC$ is a
right angle so
$\sec \frac{\pi}{4}=OC/OB$, hence $OC=\sec \frac{\pi}{4}$.
As well, $\sec \frac{\pi}{8}=OD/OC$, hence $OD=\sec \frac{\pi}{4} \cdot \sec
\frac{\pi}{8}$, etc. Then using the formula we get $OZ=BcdefgA$.}

\begin{figure}
\begin{center}
\includegraphics{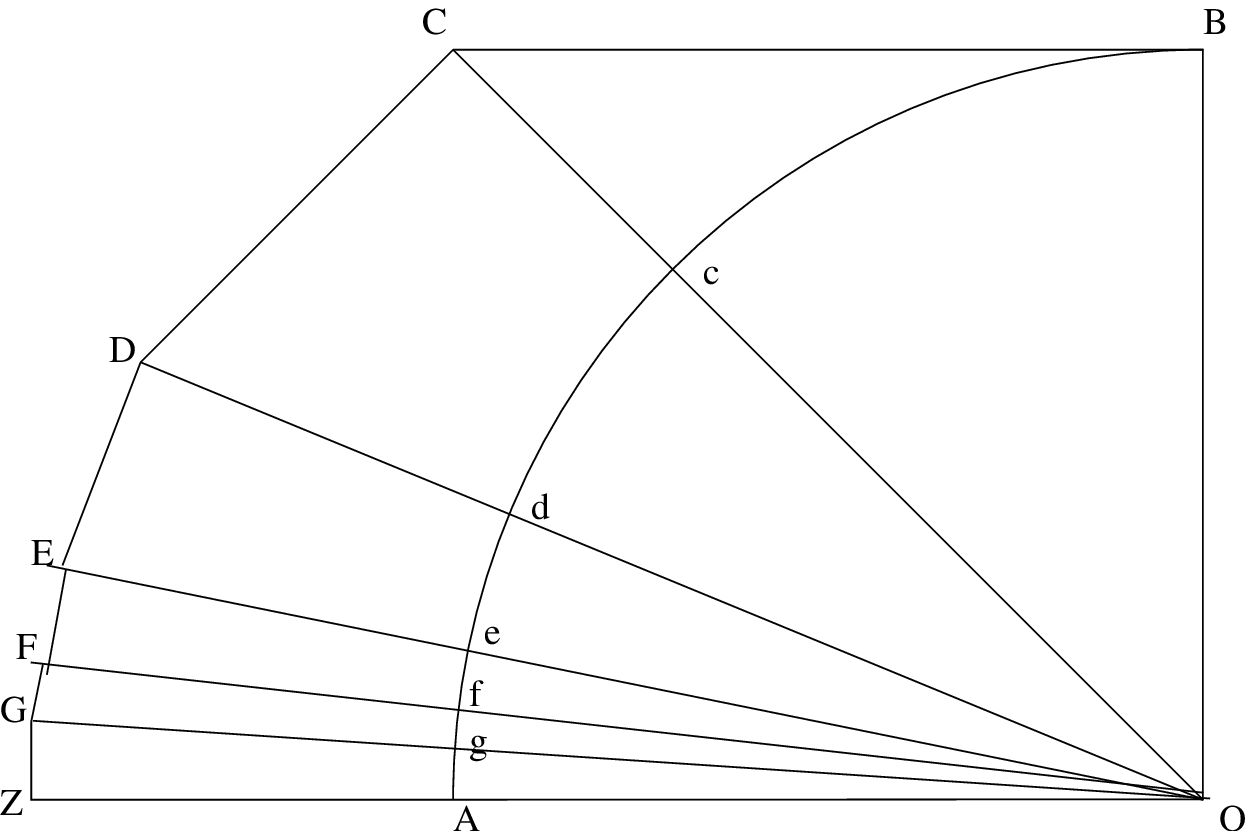}
\end{center}
\end{figure}

As well, one can easily derive many other 
constructions of this kind from the formulas of the Author.
It will be helpful to note that the points
$B,C,D,E,F$ are found on a curve such that, by putting 
any angle $AOD=\phi$ and the line $OD=v$, it turns out that $v=\frac{q \Sin \phi}{\phi}$.\footnote{Translator: This is the polar form of the quadratrix.
Concerning finding points on the quadratrix see Christoph Clavius, {\em Commentaria in Euclidis Elementa}, appendix to book VI, pp. 296-304. Available online at http://mathematics.library.nd.edu/clavius/;
also see pp. 317--318 of H. J. M. Bos, {\em On the representation of curves in Descartes' G\'eom\'etrie},
Arch. Hist. Exact Sci. \textbf{24} (1981), no. 4, 295--338.}
Then, taking any ratio between the angle $\phi$ and the right angle, whose
measure is the arc $q$, let $\phi=\frac{m}{n}q$;
it will be $v=\frac{n}{m} \Sin\frac{m}{n}q$ and the line $v$ can thus be
assigned geometrically.
Further, taking the angle $\phi$ as continually decreasing, it will
come finally to a vanishing angle, for which $\frac{\Sin \phi}{\phi}=1$,
and then the line $v$ is clearly equal to the quadrant $q$. One
can make any number of additional similar formulas.}

In Excerpts from the Manuscripts of {\em Descartes},\footnote{Translator: 
See
Ren\'e Descartes, {\em R. Des-Cartes opuscula posthuma, physica, et
mathematica}, 1701: part 6, Excerpta ex MSS. R. Des-Cartes, p. 6}
a certain geometric construction
which quickly approaches the true measure of the circle
is briefly described.
This construction,
which either {\em Descartes} himself had found, or which had been communicated
 by someone else,
especially at that time
indicates brilliantly the insightful character of its discoverer.
Those who later handled this same argument, as far as I know at least,
have not made mention of this extraordinary construction, so that it is in
danger that it disappear altogether into oblivion. This demonstration,
which is given with nothing added to it,
can be supplied without difficulty; truly not only the elegance
of this fertile construction merits more study, but 
the notable conclusions that can be derived from it would
altogether by themselves be worthy of attention.
This most beautiful construction 
is proposed thus in the words of {\em Descartes} himself:

\begin{quotation}
\marginpar{Fig. 1}
I find nothing more suitable for the quadrature of the circle than this:
if to a given square $bf$ is adjoined a rectangle $cg$
contained by the lines $ac$ and $bc$ which is equal to a fourth
part of the square $bf$: likewise a rectangle $dh$ is made from the lines
$da,dc$, equal to a fourth part of the preceding; and in the same way
a rectangle $ei$, and further infinitely many others on to $x$: and this line
$ax$
will be the diameter of a circle whose circumference is equal to the
perimeter of the square $bf$.
\end{quotation}

\begin{figure}
\begin{center}
{\Large Fig. 1}\\
\includegraphics{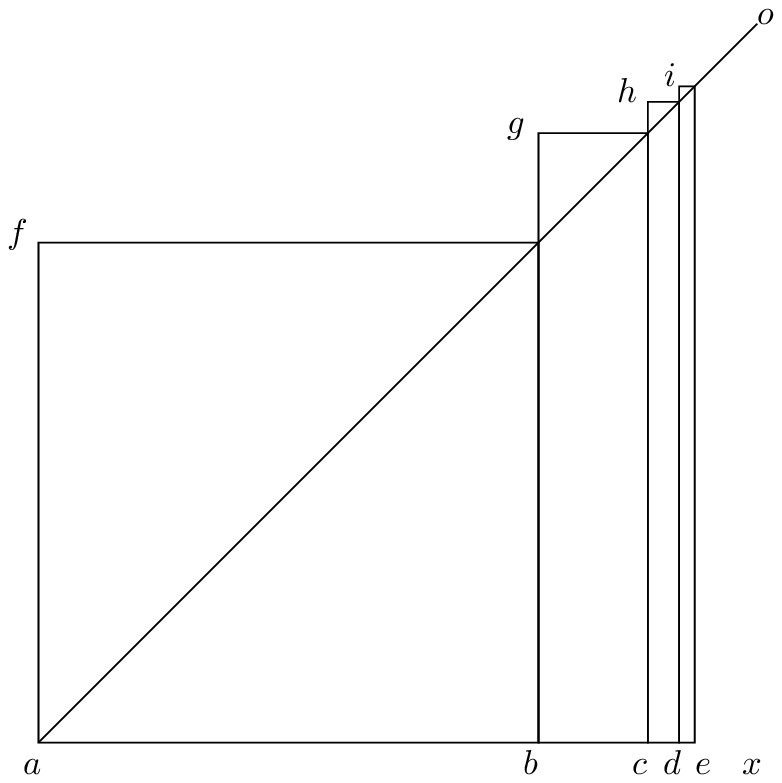}
\end{center}
\end{figure}

The strength of this method therefore consists in that by continually
adjoining rectangles of this type, $cg,dh,ei$, etc., 
whose top right angles fall on the extended diagonal of the square,
finally leading to the point $x$ at 
which ends 
the diameter $ax$ of a circle 
whose circumference is equal to the perimeter of the square $bf$, or
four times the line $ab$.

Since each of these rectangles is equal to a fourth part of the preceding one,
as is already observed by {\em Descartes} himself,
it is clear that the sum of all these rectangles will be equal to a third part
of the square $bf$; indeed this is clear because the sum of the series
\[
\frac{1}{4}+\frac{1}{16}+\frac{1}{64}+\frac{1}{256}+\etc
\]
continued to infinity
is $=\frac{1}{3}$.

{\em Descartes} further indicates the rule on which this construction
rests;
he begins namely with regular polygons of $8,16,32,64$, etc. sides,
whose perimeters are mutually equal to the perimeter of the square $bf$.
Now as $ab$ is the diameter of a circle inscribed in this square,
it is thus affirmed that $ac$ is the diameter of a circle inscribed in
an octagon,
and indeed that $ad$ and $ae$ are the diameters of circles inscribed in
a $16$gon and $32$gon respectively, and so on. Thus one sees that $ax$
is the diameter of a circle inscribed in a polygon of infinitely many sides,
whose circumference
will therefore be equal to the perimeter of the square.

For clarity, I shall elaborate the demonstration of this construction.
I observe that
what are spoken of presently as the diameters of circles can equally
well serve as radii, so that
$ab,ac,ad,ae$, etc. can be seen to be radii of circles about which if regular
polygons
of $4,8,16,32$, etc. sides are circumscribed, the perimeters of the polygons will be equal.

\begin{center}
{\Large Problem}
\end{center}

{\em Given a circle about which a regular polygon has been circumscribed,
to find another circle 
such that if a regular polygon with twice as many
sides circumscribes it, the perimeter of the first polygon will be equal
to the perimeter of the second polygon.}

\begin{center}
{\Large Solution}
\end{center}

\marginpar{Fig. 2}

Let $ENM$ be the given circle, with center at $C$, and $EP$ half of one side of
the circumscribing polygon; further let $CF$ be the radius of
the circle that is being sought, and $FQ$ half of one side of the polygon
which is to circumscribe it. It is therefore necessary that $FQ$ 
be half of $EP$ and that the angle $FCQ$ be half of the angle $ECP$.
Let the line $CQ$ bisect the angle $ECP$, and the line $QO$
parallel to $CE$ bisect the line $EP$. Now since\footnote{Translator:
Ed Sandifer in his May 2008 {\em How Euler did it}
works through the details of this paper. In particular,
Sandifer explains that one can show $EV/CE=EP(CE+CP)$ using 
the identity $\cot \frac{\theta}{2}=\cot \theta +\sec \theta$; see pp.
151--152 of I. M. Gelfand and M. Saul, {\em Trigonometry} for this identity.}

\begin{eqnarray*}
&&EV:CE=FQ:CF\\
\textrm{and}&&EV:CE=EP:CE+CP\\
\textrm{it will be}&&FQ:CF=EP:CE+CP
\end{eqnarray*}
but because $FQ=\frac{1}{2}EP$, it will further be
\[
CF=\frac{1}{2}(CE+CP).
\]
Then taking away $CF$ one will have
\[
EF=\frac{1}{2}(CP-CE)
\]
from which the
rectangle will be
\[
CF\cdot EF=\frac{1}{4}(CP^2-CE^2)=\frac{1}{4}EP^2,
\]
and so the point $F$
should be defined such that the rectangle contained
by $CF$ and $EF$ is equal to a fourth part of the square of the line $EP$,
or to the square itself of the line $FQ$.

\begin{figure}
\begin{center}
{\Large Fig. 2}\\
\includegraphics{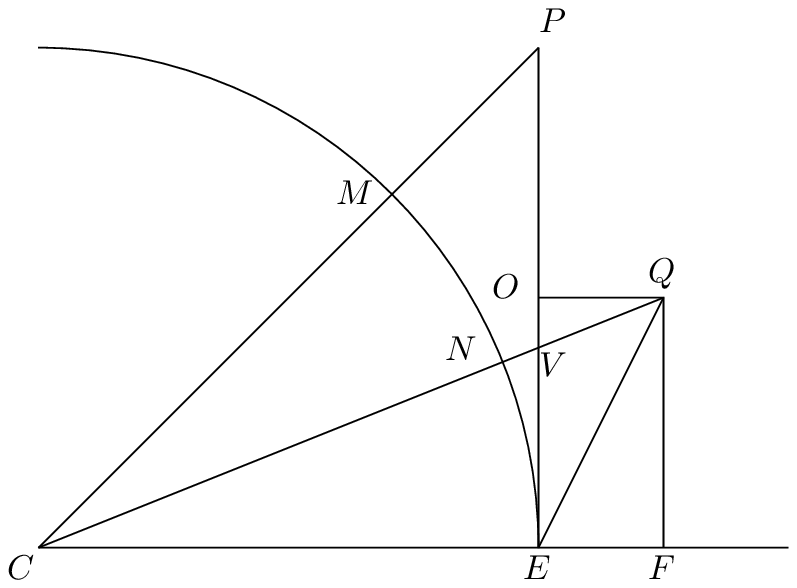}
\end{center}
\end{figure}

\begin{center}
{\Large Corollary 1}
\end{center}

Since $CF\cdot EF=FQ^2$, it will be $CF:FQ=FQ:EF$, 
whence having extended the line $QE$, the triangle $FQE$ will be similar
to the triangle $FCQ$, or $ECV$, and hence the angle $FQE$ will be equal
to the angle $ECV$.

\begin{center}
{\Large Corollary 2}
\end{center}

Since $CE:EV=EO:EF$, the point $F$ may also be defined thus:
from the point $O$ the line normal to the line $CV$ extended is drawn,
and it will meet the base $CE$ at $F$.

\begin{center}
{\Large Corollary 3}
\end{center}

If the polygon inscribed in the circle $ENM$ has $n$ sides, 
the angle $ECP$ will be $=\frac{\pi}{n}$,
where $\pi$ denotes the measure of two right angles; and the angle
$FCQ=\frac{\pi}{2n}$. Hence if the radius $CE=r$, it will be
\[
EP=r\Tan\frac{\pi}{n} \quad  \textrm{and} \quad
FQ=\frac{1}{2}r\Tan\frac{\pi}{n}.
\]

\begin{center}
{\Large Corollary 4}
\end{center}

Now because the angle $FQE=\frac{\pi}{2n}$, it will be
\[
EF=FQ\Tan\frac{\pi}{2n}=\frac{1}{2}r\Tan\frac{\pi}{n}\Tan\frac{\pi}{2n}.
\]
Indeed now, if we call
$CF=s$, it will be
\[
FQ=s\Tan\frac{\pi}{2n},
\]
whence because
\[
FQ=\frac{1}{2}r\Tan\frac{\pi}{n},
\]
it will be
\[
s=\frac{1}{2}r\Tan\frac{\pi}{n}\Cot\frac{\pi}{2n}.
\]

\begin{center}
{\Large Demonstration of Descartes' construction}
\end{center}

\marginpar{Fig. 3}
Here let $CE$ be the radius of a circle inscribed in a square,\footnote{Translator: Here we are inscribing the circle in $n$-gons which all have the same perimeter.}
$CF$ inscribed in an octagon, $CG$ in a regular polygon of $16$ sides,
$CH$ in a polygon of $32$ sides and so on. Next let $EP$ be half of one side
of
the square, $FQ$ half of one side of the octagon,
$GR$ half of one side of the polygon
with $16$ sides, $HS$ half of one side of the polygon with $32$ sides, etc.,
and because these polygons are assumed to have the same perimeter,
it will be
\[
FQ=\frac{1}{2}EP, \quad GR=\frac{1}{2}FQ=\frac{1}{4}EP,\quad
HS=\frac{1}{2}GR=\frac{1}{4}FQ=\frac{1}{8}EP, \quad \textrm{etc.}
\]
Now from the previous problem, we have $CF\cdot EF=\frac{1}{4}EP^2=FQ^2$;
indeed, we get in the same way from it
\[
\begin{split}
&CG\cdot FG=\frac{1}{4}FQ^2=\frac{1}{4}CF\cdot EF=GR^2,\\
&CH\cdot GH=\frac{1}{4}GR^2=\frac{1}{4}CG\cdot FG=HS^2 \etc
\end{split}
\] 
and the points
$F,G,H$, etc. 
are plainly determined in the same way as was done in 
{\em Descartes'} construction;
and because
the intervals $EF,FG,GH$, etc. are made continually smaller,
the final point $x$ will be approached quickly enough. Then $Cx$ will be
the radius of the circle whose circumference is equal to the perimeter of
the preceding polygons, and thus to eight times the line $EP$. \newline Q. E. D

\begin{figure}
\begin{center}
{\Large Fig. 3}\\
\includegraphics{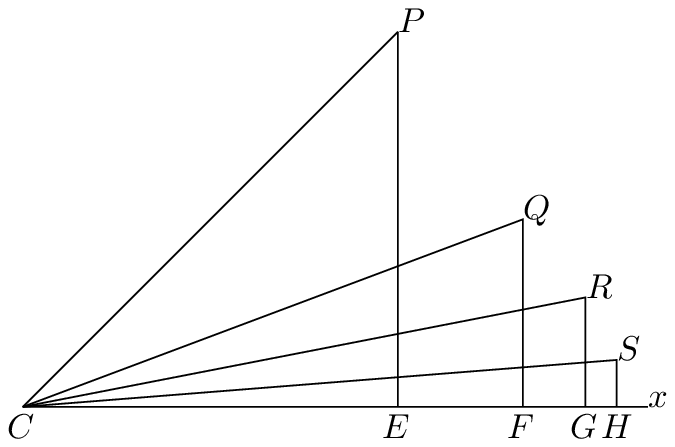}
\end{center}
\end{figure}

\begin{center}
{\Large Corollary 1}
\end{center}

If one puts $CE=a$, $CF=b$, $CG=c$, $CH=d$, etc., the progression
of these quantities is such that, as $EP=a$,
\[
b(b-a)=\frac{1}{4}aa,\quad c(c-b)=\frac{1}{4}b(b-a),\quad
d(d-c)=\frac{1}{4}c(c-b) \quad \etc
\]
and hence
\[
b=\frac{a+\surd 2aa}{2}, \quad c=\frac{b+\surd(2bb-ab)}{2},
\quad 
d=\frac{c+\surd(2cc-bc)}{2} \quad \etc
\]
and the limit of these quantities is the radius of a circle whose
circumference is $=8a$.

\begin{center}
{\Large Corollary 2}
\end{center}

Since the angle $ECP$ is half a right angle, or
\[
ECP=\frac{\pi}{4},
\]
the angles will be
\[
FCQ=\frac{\pi}{8}, \quad GCR=\frac{\pi}{16}, \quad
HCS=\frac{\pi}{32}, \quad \etc
\]
Hence because 
\[
EP=a, \quad FQ=\frac{1}{2}a, \quad GR=\frac{1}{4}a, \quad
HS=\frac{1}{8}a, \quad \etc
\]
with cotangents it will be
\[
CE=a\Cot\frac{\pi}{4}, \quad CF=\frac{1}{2}a\Cot\frac{\pi}{8}, \quad
CG=\frac{1}{4}a\Cot\frac{\pi}{16}, \quad CH=\frac{1}{8}a\Cot\frac{\pi}{32}
\quad \etc
\]
Then with $n$ denoting an infinite number, the last of these lines would
\[
=\frac{1}{n}\Cot\frac{\pi}{4n}.
\]

\begin{center}
{\Large Corollary 3}
\end{center}

But $\Cot\frac{\pi}{4n}=1:\Tan\frac{\pi}{4n}$; and because the angle
$\frac{\pi}{4n}$ is infinitely small, it will be
\[
\Tan\frac{\pi}{4n}=\frac{\pi}{4n} \quad \textrm{and hence} \quad
\Cot\frac{\pi}{4n}=\frac{4n}{\pi}.
\]
So were the last of these lines $=\frac{4a}{\pi}$, then
the circumference of a circle described with this radius will be $=2\pi \cdot \frac{4a}{\pi}=8a$.

\begin{center}
{\Large Corollary 4}
\end{center}

Then, since by corollary 4 of the preceding problem\footnote{Translator: I don't see quite what corollary 4 of the preceding problem has to do with this. We
know that $CF\cdot EF=FQ^2, CG\cdot FG=GR^2$, etc. from the Demonstration,
from which $EF=FQ\Tan FCQ, FG=GR \Tan GCR$, etc. follow by similar triangles.}
\[
EF=FQ\Tan FCQ,
\]
by the same rule it will be
\[
FG=GR \Tan GCR, \quad GH=HS\Tan HCS \quad \etc
\]
from which the other intervals can be expressed in the following way
\[
EF=\frac{1}{2}a\Tan\frac{\pi}{8}, \quad FG=\frac{1}{4}a\Tan\frac{\pi}{16},
\quad GH=\frac{1}{8}a\Tan\frac{\pi}{32} \quad \etc,
\]
with the first indeed in analogy
\[
CE=a\Tan\frac{\pi}{4}=a.
\]

\begin{center}
{\Large Corollary 5}
\end{center}

With these combined with the preceding we will find
\[
\begin{split}
&CF=a\Big(\Tan\frac{\pi}{4}+\frac{1}{2}\Tan\frac{\pi}{8}\Big)=\frac{1}{2}a\Cot\frac{\pi}{8},\\
&CG=a\Big(\Tan\frac{\pi}{4}+\frac{1}{2}\Tan\frac{\pi}{8}+\frac{1}{4}\Tan\frac{\pi}{16}\Big)=\frac{1}{4}a\Cot\frac{\pi}{16},\\
&CH=a\Big(\Tan\frac{\pi}{4}+\frac{1}{2}\Tan\frac{\pi}{8}+
\frac{1}{4}\Tan\frac{\pi}{16}+\frac{1}{8}\Tan\frac{\pi}{32}\Big)=\frac{1}{8}a\Cot\frac{\pi}{32}\\
&\etc
\end{split}
\]
and so in this way, the sums of all these progressions can be
assigned without any difficulty.

\begin{center}
{\Large Corollary 6}
\end{center}

Therefore progressing to infinity we will obtain the summation of the series
\[
\Tan\frac{\pi}{4}+\frac{1}{2}\Tan\frac{\pi}{8}+\frac{1}{4}\Tan\frac{\pi}{16}+\frac{1}{8}\Tan\frac{\pi}{32}+\etc=\frac{4}{\pi},
\]
which is thus obtained by the quadrature of the circle. I shall take occasion
to solve the following problem.

\begin{center}
{\Large Problem}
\end{center}

{\em With $\phi$ denoting the arc of a circle whose radius is $=1$, to find the
sum of the infinite series}
\[
\Tan \phi+\frac{1}{2}\Tan \frac{1}{2}\phi + \frac{1}{4}\Tan \frac{1}{4}\phi+
\frac{1}{8}\Tan \frac{1}{8}\phi+\frac{1}{16}\Tan\frac{1}{16}\phi+\etc
\]

\begin{center}
{\Large Solution}
\end{center}

If in Fig. 2, as constructed above, one sets
the angle $ECP=\phi$, it will be $FCQ=\frac{1}{2}\phi$:
now putting $FQ=1$ it will be $EP=2$, and hence
\[
CE=2\Cot\phi, \quad CF=\Cot\frac{1}{2}\phi \quad \textrm{and} \quad
EF=\Tan\frac{1}{2}\phi,
\]
from
which one has
\[
2\Cot\phi=\Cot\frac{1}{2}\phi-\Tan\frac{1}{2}\phi
\quad \textrm{and} \quad \Tan\frac{1}{2}\phi=\Cot\frac{1}{2}\phi-2\Cot\phi.
\]
and in the same way
\[
\Tan\phi=\Cot\phi-2\Cot 2\phi.
\]
The values of the tangents found in the given series can be written
by cotangents as
\begin{eqnarray*}
\Tan\phi&=&\Cot\phi-2\Cot 2\phi,\\
\frac{1}{2}\Tan\frac{1}{2}\phi&=&\frac{1}{2}\Cot\frac{1}{2}\phi-\Cot\phi,\\
\frac{1}{4}\Tan\frac{1}{4}\phi&=&\frac{1}{4}\Cot\frac{1}{4}\phi-\frac{1}{2}\Cot\frac{1}{2}\phi,\\
\frac{1}{8}\Tan\frac{1}{8}\phi&=&\frac{1}{8}\Cot\frac{1}{8}\phi-\frac{1}{4}\Cot\frac{1}{4}\phi\\
&&\etc
\end{eqnarray*}
and adding them together we will get
\begin{eqnarray*}
\Tan\phi=\Cot\phi-2\Cot 2\phi,\\
\Tan\phi+\frac{1}{2}\Tan\frac{1}{2}\phi=\frac{1}{2}\Cot\frac{1}{2}\phi-2\Cot 2\phi,\\
\Tan\phi+\frac{1}{2}\Tan\frac{1}{2}\phi+\frac{1}{4}\Tan\frac{1}{4}\phi=\frac{1}{4}\Cot\frac{1}{4}\phi-2\Cot 2\phi,\\
\Tan\phi+\frac{1}{2}\Tan\frac{1}{2}\phi+\frac{1}{4}\Tan\frac{1}{4}\phi+\frac{1}{8}\Tan\frac{1}{8}\phi=\frac{1}{8}\Cot\frac{1}{8}\phi-2\Cot 2\phi\\
\etc,
\end{eqnarray*}
which continued to infinity will be
$\frac{1}{n}\Cot\frac{1}{n}\phi=\frac{1}{\phi}$ if $n$ denotes an infinity
number, because $\Tan \frac{1}{n}\phi=\frac{\phi}{n}$ and hence
$\Cot \frac{1}{n}\phi=\frac{n}{\phi}$. Hence the sum of the given series
will be
\[
\Tan\phi+\frac{1}{2}\Tan\frac{1}{2}\phi+\frac{1}{4}\Tan\frac{1}{4}\phi+\frac{1}{8}\Tan\frac{1}{8}\phi+\etc=\frac{1}{\phi}-2\Cot 2\phi.
\]
Whence if $2\phi$ is a right angle or $\phi=\frac{\pi}{4}$,
then since $\Cot\frac{\pi}{2}=0$, the sum of the series would be $=\frac{1}{\phi}=\frac{4}{\pi}$,
which was treated in the above case.

From this series many others can be derived which are no less noteworthy.

I. By differentiating this series we obtain
\[
\frac{1}{\Cos\phi^2}+\frac{1}{4\Cos\frac{1}{2}\phi^2}+
\frac{1}{4^2\Cos\frac{1}{4}\phi^2}+\frac{1}{4^3\Cos\frac{1}{8}\phi^2}+
\frac{1}{4^4\Cos\frac{1}{16}\phi^2}+\etc=-\frac{1}{\phi\phi}+\frac{4}{\Sin 2\phi^2};
\]
or, since $\frac{1}{\Cos\phi}=\Sec\phi$ it will also be
\[
(\Sec\phi)^2+\frac{1}{4}\Big(\Sec\frac{1}{2}\phi\Big)^2
+\frac{1}{4^2}\Big(\Sec\frac{1}{4}\phi\Big)^2+\frac{1}{4^3}\Big(\Sec\frac{1}{8}\phi\Big)^2+\etc=\frac{1}{\Sin\phi^2 \Cos\phi^2}-\frac{1}{\phi\phi}
\]

II. Next, because $\Cos\phi^2=\frac{1+\Cos 2\phi}{2}$ and $\Sin 2\phi^2=
\frac{1-\Cos 4\phi}{2}$, by dividing everything by two
it will be
\[
\begin{split}
&\frac{1}{1+\Cos 2\phi}+\frac{1}{4(1+\Cos\phi)}+
\frac{1}{4^2(1+\Cos\frac{1}{2}\phi)}+\frac{1}{4^3(1+\Cos\frac{1}{4}\phi)}+\etc\\
&=\frac{2}{1-\Cos 4\phi}-\frac{1}{2\phi\phi}
\end{split}
\]
or by writing $\frac{1}{2}\phi$ for $\phi$
\[
\begin{split}
&\frac{1}{1+\Cos\phi}+\frac{1}{4(1+\Cos\frac{1}{2}\phi)}+
\frac{1}{4^2(1+\Cos\frac{1}{4}\phi)}+\frac{1}{4^3(1+\Cos\frac{1}{8}\phi)}+\etc\\
&=\frac{2}{1-\Cos 2\phi}-\frac{2}{\phi\phi}.
\end{split}
\]

III.  If the series which has been found is multiplied by $d\phi$ and
integrated, because
\[
\int d\phi \Tan\phi=\int \frac{d\phi \Sin\phi}{\Cos\phi}=
-l\Cos\phi \quad \textrm{and} \quad \int 2d\phi\Cot 2\phi=l\Sin 2\phi,
\]
one will have
\[
-l\Cos\phi-l\Cos\frac{1}{2}\phi-l\Cos\frac{1}{4}\phi-l\Cos\frac{1}{8}\phi-l\Cos\frac{1}{16}\phi-\etc=l\phi-l\Sin 2\phi+\textrm{Const.}
\]
In order to define this constant, let us put $\phi=0$, and because
$l\Cos 0=l1=0$ we will have $0$ for the first part, while for the
second part, since $\Sin 2\phi=2\phi$, we will have $l\phi-l2\phi+\textrm{Const.}=-l2+\textrm{Const.}$, whence $\textrm{Const.}=l2$. 
Hence switching to numbers instead of the logarithms of numbers it will be
\[
\frac{1}{\Cos\phi\Cos\frac{1}{2}\phi\Cos\frac{1}{4}\phi\Cos\frac{1}{8}\phi\Cos\frac{1}{16}\Phi\etc}=\frac{2\phi}{\Sin 2\phi}
\]

IV. Since
\[
\frac{1}{\Cos\phi}=\Sec\phi,
\]
this theorem can also be expressed according to secants as
\[
\Sec\phi \, \Sec\frac{1}{2}\phi \, \Sec\frac{1}{4}\phi \, \Sec\frac{1}{8}\phi
\, \Sec\frac{1}{16}\phi\,\etc=\frac{2\phi}{\Sin 2\phi}.
\]
From this, if the ratio of the diameter to the circumference is put $=1:\pi$
and $q$ denotes a right angle, if we set $2\phi=q=\frac{\pi}{2}$ it will
be
\[
\Sec\frac{1}{2}q \, \Sec\frac{1}{4}q \, \Sec\frac{1}{8}q \, \Sec\frac{1}{16}q
\, \Sec\frac{1}{32}q \, \etc=\frac{\pi}{2}.
\]

\begin{center}
{\Large Problem}
\end{center}

{\em To find a series of quantities: $a,b,c,d,e,f$, etc.
which have the property that}
\[
c(c-b)=\frac{1}{4}b(b-a), \quad d(d-c)=\frac{1}{4}c(c-b), \quad e(e-d)=\frac{1}{4}d(d-c) \quad \etc
\]
{\em or such that the quantities thence formed}
\[
b(b-a),\quad c(c-b), \quad d(d-c),\quad e(e-d),\quad f(f-e),\quad \etc
\]
{\em decrease in quadruple ratio.}

\begin{center}
{\Large Solution}
\end{center}

Since
$\Tan\frac{1}{2}\phi=\Cot\frac{1}{2}\phi-2\Cot\phi$,\footnote{Translator:
This follows from the addition formula for $\tan$: $\tan(\alpha+\beta)=
\frac{\tan\alpha+\tan\beta}{1-\tan\alpha \tan\beta}$.} if
we multiply each side by
$\Cot\frac{1}{2}\phi$, since $\Tan\frac{1}{2}\phi\Cot\frac{1}{2}\phi=1$
it will be
\[
\Cot\frac{1}{2}\phi\Big(\Cot\frac{1}{2}\phi-2\Cot\phi\Big)=1.
\]
Let us then set
\[
a=r\Cot\phi, \quad b=\frac{1}{2}r\Cot\frac{1}{2}\phi,\quad
c=\frac{1}{4}r\Cot\frac{1}{4}\phi,\quad
d=\frac{1}{8}r\Cot\frac{1}{8}\phi\quad \etc,
\]
and it will be
\begin{eqnarray*}
\frac{2b}{r}\Big(\frac{2b}{r}-\frac{2a}{r}\Big)=1,&\textrm{hence}&b(b-a)=\frac{rr}{4},\\
\frac{4c}{r}\Big(\frac{4c}{r}-\frac{4b}{r}\Big)=1,&\textrm{hence}&c(c-b)=\frac{rr}{4^2},\\
\frac{8d}{r}\Big(\frac{8d}{r}-\frac{8c}{r}\Big)=1,&\textrm{hence}&d(d-c)=\frac{rr}{r^3}\\
&\etc&
\end{eqnarray*}
whence this series
\[
a=r\cot\phi, \quad b=\frac{1}{2}r\Cot\frac{1}{2}\phi, \quad
c=\frac{1}{4}r\Cot\frac{1}{4}\phi, \quad
d=\frac{1}{8}r\Cot\frac{1}{8}\phi \quad \etc
\]
has the property that the quantities thence formed 
\[
b(b-a), \quad c(c-b), \quad d(d-c), \quad e(e-d) \quad \etc
\]
decrease in quadruple ratio.

\begin{center}
{\Large Corollary 1}
\end{center}

Given the first two terms $a$ and $b$, all the remaining $c,d,e,f$
are thence successively determined, such that
\[
c=\frac{b+\surd(2bb-ab)}{2}, \quad
d=\frac{c+\surd(2cc-bc)}{2}, \quad e=\frac{d+\surd(2dd-cd)}{2} \quad \etc
\]
and hence with the first two terms taken at our pleasure,
the entire series can be exhibited by means of these formulas.

\begin{center}
{\Large Corollary 2}
\end{center}

Moreover, with the terms $a$ and $b$ given, the angle $\phi$ and the
quantity $r$
can be thus defined from them
\[
\Tan\phi=\frac{2\surd(bb-ab)}{a} \quad \textrm{and} \quad r=2\surd(bb-ab);
\]
then, having found the angle $\phi$ all the remaining terms can also
be expressed, as
\[
c=\frac{1}{4}r\Cot\frac{1}{4}\phi,
\quad d=\frac{1}{8}r\Cot\frac{1}{8}\phi,
\quad e=\frac{1}{16}r\Cot\frac{1}{16}\phi \quad \etc
\]

\begin{center}
{\Large Corollary 3}
\end{center}

Hence the infinitesimal terms of this series will be
$=\frac{r}{\phi}$,
to which value the terms of the series converge fairly
quickly.\footnote{Translator: $\cot x=\frac{1}{x}-\frac{x}{3}-\frac{x^3}{45}
-\cdots$, so the difference between any of $c,d,e$, etc. and $\frac{r}{\phi}$
is $O(\phi)$.}
Namely an arc is sought in the circle with radius $=1$, whose
tangent
\[
=\frac{2\surd(bb-ab)}{a};
\]
let this arc be $=\phi$, and the infinitesimal terms of our series
will be
\[
=\frac{2\surd(bb-ab)}{\phi}.
\]

\begin{center}
{\Large Scholion}
\end{center}

\marginpar{Fig. 3}

It will be useful to note here that the points $P,Q,R,S,x$
(Fig. 3)
are situated on the quadratrix of antiquity,
because the 
line segments
$EP,FQ,GR,HS$
have the same ratio to each other
as the angles $ECP,FCQ,GCR,HCS$, etc.\footnote{Translator: Say
$x_1=y_1\cot(\frac{y_1 \pi}{2})$, i.e. $(x_1,y_1)$ is on the quadratrix.
Using polar coordinates,
$\frac{x_1}{y_1}=\cot \theta_1$. Suppose that
$y_2=\frac{y_1}{2}$ and $\theta_2=\frac{\theta_1}{2}$. Writing $(x_2,y_2)$ in polar
coordinates we get $\frac{x_2}{y_2}=\cot \theta_2$. Then, in a few lines one
can show
using $x_1=y_1\cot(\frac{y_1 \pi}{2})$ and $\frac{x_1}{y_1}=\cot \theta_1$
that $x_2=y_2\cot(\frac{y_2 \pi}{2})$, i.e. that $(x_2,y_2)$ is on the quadratrix.}
And since $x$, which is where this curve intersects the
base,
here as before has been found to indicate the quadrature of the circle,
which is the very reason
for the name of this curve,
the construction of {\em Descartes} 
agrees singularly with this quadrature of antiquity,
but it offers more conveniently and accurately the
points $E,F,G,H$, etc. in succession,
than what one could hope for by the continual bisection of
angles.\footnote{Translator: I don't see what Euler was thinking about when he
says that Descartes' method is better than the one based on continued bisection
of the angle. This statement needs clarification. References are
given in the footnotes of the {\em Opera omnia}, I.15, pp. 1--15.}

\end{document}